%% file: articulo1.tex
\documentclass[12pt]{amsart}

\usepackage[latin1]{inputenc}
\usepackage[active]{srcltx}

\usepackage[all]{xy}

\usepackage{cases}

\usepackage{enumerate} \usepackage{float} \textwidth=6truein
\usepackage{amssymb}
\newcommand{\Z}{{\mathbb Z}} \newcommand{\Q}{{\mathbb Q}}

\newcommand{\liminv}{\operatornamewithlimits{\hbox{$\lim$}}}
\newcommand{\limdir}{\operatornamewithlimits{\hbox{$\colim$}}}
\newcommand{\colim}{\operatornamewithlimits{colim}}

\newcommand{\Hom}{\operatorname{Hom}\nolimits}

\renewcommand{\ker}{\operatorname{Ker}\nolimits}
\newcommand{\im}{\operatorname{Im}\nolimits}
\newcommand{\coker}{\operatorname{Coker}\nolimits}
\newcommand{\coim}{\operatorname{Coim}\nolimits}

\newcommand{\Ab}{\operatorname{Ab}\nolimits}

\newcommand{\Ob}{\operatorname{Ob}\nolimits}

\newcommand{\A}{\ifmmode{\mathcal{A}}\else${\mathcal{A}}$\fi}
\newcommand{\B}{\ifmmode{\mathcal{B}}\else${\mathcal{B}}$\fi}
\newcommand{\C}{\ifmmode{\mathcal{C}}\else${\mathcal{C}}$\fi}
\newcommand{\D}{\ifmmode{\mathcal{D}}\else${\mathcal{D}}$\fi}
\newcommand{\G}{\ifmmode{\mathcal{G}}\else${\mathcal{G}}$\fi}
\newcommand{\I}{\ifmmode{\mathcal{I}}\else${\mathcal{I}}$\fi}
\newcommand{\J}{\ifmmode{\mathcal{J}}\else${\mathcal{J}}$\fi}
\newcommand{\K}{\ifmmode{\mathcal{K}}\else${\mathcal{K}}$\fi}
\renewcommand{\O}{\ifmmode{\mathcal{O}}\else${\mathcal{O}}$\fi}
\renewcommand{\P}{\ifmmode{\mathcal{P}}\else${\mathcal{P}}$\fi}
\newcommand{\U}{\ifmmode{\mathcal{U}}\else${\mathcal{U}}$\fi}
\newcommand{\M}{\ifmmode{\mathcal{M}}\else${\mathcal{M}}$\fi}
\newcommand{\N}{\ifmmode{\mathcal{N}}\else${\mathcal{N}}$\fi}
\newcommand{\Ss}{\ifmmode{\mathcal{S}}\else${\mathcal{S}}$\fi}
\newcommand{\T}{\ifmmode{\mathcal{T}}\else${\mathcal{T}}$\fi}
\newcommand{\Ff}{\ifmmode{\mathcal{F}}\else${\mathcal{F}}$\fi}
\newcommand{\Ll}{\ifmmode{\mathcal{L}}\else${\mathcal{L}}$\fi}
\newtheorem{Thm}{Theorem}[section]
\newtheorem{Prop}[Thm]{Proposition}

\newtheorem{Lem}[Thm]{Lemma}
\newtheorem{Claim}{Claim}[Thm]

\newtheorem*{ThmA}{Theorem A}
\newtheorem*{ThmB}{Theorem B}
\newtheorem*{PropC}{Proposition C}

\theoremstyle{definition}
\newtheorem{Defi}[Thm]{Definition}
\newtheorem{Rmk}[Thm]{Remark}
\newtheorem{Ex}[Thm]{Example}

\theoremstyle{remark}

\theoremstyle{plain}

\newcommand{\definicio}{\stackrel{\text{def}}{=}}

\title{A family of acyclic functors}

\author{Antonio D{\'\i}az Ramos}
\address{Department of Mathematical Sciences\\ King's College\\
University of Aberdeen \\ ABERDEEN AB24 3UE U.K.}
\email{a.diaz@maths.abdn.ac.uk}
\date{\today}

\begin{document}

\maketitle

\input{introduction}


\input{graded_posets}

\input{spectral}

\input{projective}

\input{limacyI}



\input{dual_results}





\input{biblio}
\end{document}

%% file: introduction.tex
\section{Introduction and summary}\label{section_introduction}
In this paper we focus on the following problem
\begin{equation}\label{equ_intro_main_problem}
\text{find $\limdir$-acyclic objects in $\Ab^\C$.}
\end{equation}
Here, $\Ab$ denote the category of abelian groups and $\Ab^\C$
denote the (abelian) functor category for the small category $\C$.
The functor $\limdir:\Ab^\C\rightarrow \Ab$ is the direct limit
functor and $F\in \Ab^\C$ is $\limdir$-acyclic if $\limdir_i F=0$
for $i\geq 1$ (see \cite{weibel} and \cite{cohn}, and the classical
books of Cartan and Eilenberg \cite{cartan} and of MacLane
\cite{maclanehom}). It is clear that if $F$ is projective then it is
$\limdir$-acyclic but, in the same way as not every flat module is
projective (see, for example, \cite[Section 3.2]{weibel}), we may be
missing $\limdir$-acyclic objects if we just consider projective
ones.

We shall assume the hypothesis that the category $\C$ is a graded
partially ordered set (a graded poset for short). These are special
posets in which we can assign an integer to each object (called the
degree of the object) in such a way that preceding elements are
assigned integers which differs in $1$. Thus a graded poset can be
divided into a set of ``layers" (the objects of a fixed degree), and
these layers are linearly ordered. Any simplicial complex (viewed as
the poset of its simplices with the inclusions among them) and any
subdivision category is a graded poset. Also, every $CW$-complex is
(strong) homotopy equivalent to a simplicial complex, and thus to a
graded poset.

To attack problem (\ref{equ_intro_main_problem}) we start giving a
characterization of the projective objects in $\Ab^\C$. Recall that
for any small category $\C$ (not necessarily a poset), the
projective objects in $\Ab^\C$ are well known to be, by the Yoneda
Lemma, summands of direct sums of representable functors. Moreover,
if $\C$ is a poset with the descending chain condition (not
necessarily graded) then \cite[Corollary 3]{posetdim1} the
projective objects in $\Ab^\C$ are also direct sums of representable
functors (see also \cite[Proposition 7]{hoppner} and \cite[Theorem
9]{posetdim1tails} for related results). In case $\C$ is a graded
poset we characterize the projective functors in $\Ab^\C$ as those
functors which satisfy two conditions:

\begin{ThmA}[Theorem
\ref{pro_projective_bounded_graded}] Let $\C$ be a bounded below
graded poset and let $F:\C\rightarrow \Ab$ be a functor. Then $F$ is
projective if and only if:
\begin{enumerate}
\item for any object $i_0$ of $\C$ $\coker_F(i_0)$ is a free abelian group.
\item $F$ is pseudo-projective.
\end{enumerate}
\end{ThmA}
Here, $\coker_F(i_0)$ is the quotient of $F(i_0)$ by the images of
all the non-trivial morphisms arriving to $i_0$. For the actual
definition of pseudo-projectiveness see Definition
\ref{property_pseudo-projective}. The boundedness condition in the
theorem is related to the descending chain condition in the
aforementioned result, and neither of these conditions can be
dropped: consider the poset $\C=\Z$ of the integers. This graded
poset does not satisfies the descending chain condition, and thus
neither is it bounded below. The constant functor of value $\Z$ over
this poset is projective but it is not a sum of representable
functors. The constant functor of value $\Z/n$ (for some $n\geq 1$)
satisfies both conditions in Theorem A but it is not projective (see
Remark \ref{remark:neccesary_bounded}).

Theorem A is the first step towards finding $\limdir$-acyclic
objects in $\Ab^\C$. The reason is that the second of the conditions
in the theorem, i.e., pseudo-projectiveness, implies
$\limdir$-acyclicity:

\begin{ThmB}[Theorem \ref{pro_acyclic_graded}]
Let $F:\C\rightarrow \Ab$ be a pseudo-projective functor over a
bounded below graded poset $\C$. Then $F$ is $\limdir$-acyclic.
\end{ThmB}

This is the main result of this work, and it gives a family of
functors in $\Ab^\C$ which are $\colim$-acyclic but not necessarily
projective. To show that there exist functors in this situation
consider the functor
$$
\xymatrix{\Z&\Z\ar_{\times2}[l]\ar^{\times2}[r]&\Z }.
$$
This is a pseudo-projective functor which, by Theorem B, is acyclic.
Moreover, it does not satisfy condition (1) in Theorem A and so it
is not projective (see Examples \ref{examples_projective} and
\ref{limacyI_expushout}). On the other hand, pseudo-projective
functors do not cover all $\colim$-acyclic functors: the functor
$$
\xymatrix{0&\Z\ar_{0}[l]\ar^{1}[r]&\Z }.
$$
is not pseudo-projective but, as a straightforward computation
shows, it is acyclic. For vector spaces the notion of
pseudo-projectiveness becomes identical to projectiveness as
condition (1) in Theorem A is unnecessary in the context of functors
to $k-mod$ (where $k$ is a field). Even in this favorable case the
functor $$ \xymatrix{0&k\ar_{0}[l]\ar^{1}[r]&k }$$ shows that there
are acyclic functors which are not projective.

The main ingredient in the proof of Theorem B is a meticulous use of
a spectral sequence built upon the grading of the partially ordered
set $\C$:

\begin{PropC}[Proposition \ref{spectral_sequence_1}]
For a (decreasing) graded poset $\C$ and a functor $F:\C\rightarrow
\Ab$:
\begin{itemize}
\item There exists a cohomological type spectral sequence
$E^{*,*}_*$ with target $\limdir_* F$. \item There exists a
homological type spectral sequence $(E^p)_{*,*}^*$ with target the
column $E^{p,*}_1$ for each $p$.
\end{itemize}
\end{PropC}

Applications of these results to computation of integral cohomology
of posets are given in \cite{diaz}. The work also contains the dual
version of the above, in which we consider injective objects in
$\Ab^\C$, the right derived functors of the inverse limit functor
$\liminv:\Ab^\C\rightarrow \Ab$ and the respective $\liminv$-acyclic
objects.

The paper is structured as follows: in Section \ref{section
graded_posets} we introduce preliminaries about graded partially
ordered sets. In Section \ref{spectral} we build some spectral
sequences arising from the grading of a graded poset. Afterwards, in
Section \ref{section_projective}, we work out the characterization
of projective objects in $\Ab^\C$. In Section
\ref{section_Pseudo-projectivity} we prove that
pseudo-projectiveness implies $\limdir$-acyclicity. We finish with
Section \ref{section_dual_results}, where the dual definitions and
results for $\liminv$-acyclicity are stated without proof.

\textbf{Acknowledgements:} I would to thank my Ph.D. supervisor
Prof. A. Viruel for his support during the development of this work.
Also, thanks to Prof. C.A. Weibel for all his fruitful suggestions
and comments, in particular for a short proof of Lemma
\ref{lem_proj_pseudo_graded_coker_projective}.

%% file: graded_posets.tex
\section{Graded posets}\label{section graded_posets} In this section we
define a special kind of categories: graded partially ordered sets
(graded posets for short). We shall think of a poset $\P$ as a
category in which there is an arrow $p\rightarrow p'$ if and only if
$p\leq p'$. The notion of graded poset is not new and it was already
used in \cite[pp. 29-33]{dayton-weibel}. The definition there is
weaker than the one given here, being the difference that here we
ask for every morphism to factor trough morphisms of degree $1$
(some kind of ``saturation'' condition).

\begin{Defi}
If $\P$ is a poset and $p<p'$ then $p$ \emph{precedes} $p'$ if
$p\leq p''\leq p'$ implies that $p=p''$ or $p'=p''$.
\end{Defi}
\begin{Defi}
Let $\P$ be a poset. $\P$ is called \emph{graded} if there is a
function $deg:\Ob(\P)\rightarrow \Z$, called the \emph{degree
function} of $\P$, which is order preserving and that satisfies that
if $p$ precedes $p'$ then $deg(p')=deg(p)+1$. If $p$ is an object of
$\P$ then $deg(p)$ is called the \emph{degree} of $p$.
\end{Defi}

Notice that the degree function associated to a graded poset is not
unique (consider the translations $deg'=deg+c_k$ for $k\in \Z$).
According to the definition the degree function increases in the
direction of the arrows: we say that this degree function is
\emph{increasing}. If the degree function is order reversing and
satisfies the alternative condition that $p$ precedes $p'$ implies
$deg(p')=deg(p)-1$, i.e., $deg$ decreases in the direction of the
arrows, then we say that $deg$ is a \emph{decreasing} degree
function. Clearly both definitions are equivalent (by taking
$deg'=-deg$).

\begin{Ex}
The ``pushout category" $b\leftarrow a\rightarrow c$, the
``telescope category" $a\rightarrow b\rightarrow c\rightarrow...$,
and the opposite ``telescope category" $...\rightarrow c\rightarrow
b\rightarrow a$ are graded posets. The integers $\Z$ is a graded
poset. The rationals $\Q$ with the usual order is a poset but it is
not a graded poset.
\end{Ex}

If $\P$ is a graded poset and $p<p'$ then it is straightforward that
the number $deg(p')-deg(p)$ does not depend on the degree function
$deg$. Thus, we can ``extend" the degree function $deg$ to the
morphisms set $\Hom(\P)$ by $deg(p\rightarrow p')=|deg(p')-deg(p)|$.
Whenever $\P$ is a graded poset we denote by $\Ob_n(\P)$ the objects
of degree $n$ and by $\Hom_n(\P)$ the arrows of degree $n$. 

\subsection{Boundedness on graded posets.}\label{graded_poset_boundedness}
Often we will restrict to:

\begin{Defi}
A graded poset $\P$ with increasing degree function $deg$ is
\emph{bounded below} (\emph{bounded above}) if the set
$deg(\P)\subset \Z$ has a lower bound (an upper bound).
\end{Defi}

If the degree function  $deg$ of $\P$ is decreasing then $\P$ is
\emph{bounded below} (\emph{bounded above}) if and only if
$deg(\P)\subset \Z$ has an upper bound (a lower bound). If $\P$ is
bounded below and over then $N\definicio max(deg(\P))-min(deg(\P))$
exists and it is finite, and it does not depend on the degree
function $deg$. We call it the dimension of $\P$, and we say that
$\P$ is $N$ dimensional.

\begin{Ex}
The ``pushout category" $b\leftarrow a\rightarrow c$ is
$1$-dimensional, the ``telescope category" $a\rightarrow
b\rightarrow c\rightarrow ..$ is bounded below but it is not
bounded over. The opposite ``telescope category" $..\rightarrow
c\rightarrow b\rightarrow a$ is bounded over but it is not bounded
below.
\end{Ex}

Notice that in a bounded above (below) graded poset there are
maximal (minimal) elements, but that the existence of maximal
(minimal) objects does not guarantee boundedness. Also it is clear
that, in general, neither dcc posets are graded nor graded posets
are dcc.

%% file: spectral.tex
\section{A spectral sequence}\label{spectral} In this
section we shall construct spectral sequences with targets
$\limdir_i F$ and $\liminv^i F$ for $F:\C\rightarrow \Ab$ with $\C$
a graded poset. Some conditions for (weak) convergence shall be
given. We build the spectral sequences starting from filtered
differential modules (see \cite{mccleary}, where the notion of weak
convergence we use is also given).

Recall that (see \cite[Appendix II.3]{calculus}, \cite[XII.5.5]{BK},
\cite[XI.6.2]{BK} or \cite[p.409ff.]{goerss}) there is, for any
small category $\C$ and covariant functor $F:\C\rightarrow \Ab$, a
concrete chain (cochain) complex $C_*(\C,F)$ ($C^*(\C,F)$) whose
homology groups (cohomology groups) are precisely the left derived
functors $\limdir_i$ (right derived functors $\liminv^i$).

Let $N\C$ denote the nerve of the small category $\C$ whose
$n$-simplices are chain of composable morphisms in $\C$:
$\sigma=\xymatrix{ \sigma_0 \ar[r]^{\alpha_1} & \sigma_1
\ar[r]^{\alpha_2}&...\ar[r]^{\alpha_{n-1}}&\sigma_{n-1}\ar[r]^{\alpha_n}
& \sigma_n}$. Then
$$C_n(\C,F)=\bigoplus_{\sigma \in {N\C}_n} F_\sigma,$$
where $F_{\sigma}=F(\sigma_o)$. Moreover, $C_*(\C,F)$ is a
simplicial abelian group with face and degeneracy maps induced by
those of the nerve $N\C$. The chain complex $(C_*(\C,F),d)$ with
differential of degree $-1$ $d=\sum_{i=0}^n (-1)^i d_i$ satisfies
\begin{equation}
\label{Defi:limi} {\limdir}_i F= H_i(C_*(\C,F),d).
\end{equation}

For the inverse limit $\liminv:\Ab^\C\rightarrow \Ab$ there is a
cosimplicial abelian group with simplices
$$    C^n(\C,F)=\prod_{\sigma \in {N\C}_n} F^\sigma,$$
where $F^{\sigma}=F(\sigma_n)$. This cosimplicial object gives rise
to a cochain complex $(C^*(\C,F),d)$ with differential of degree $1$
$d=\sum_{i=0}^{n+1} (-1)^i d^i$. It is well known that
\begin{equation} \label{Defi:liminvi}
{\liminv}^i F= H^i(C^*(\C,F),d).
\end{equation}
%
\begin{Rmk}\label{simp_cosimp_normalization}
We can apply the Dold-Kan correspondence (see \cite[8.4]{weibel}) to
the simplical and cosimplicial abelian groups constructed above.
This means that we shall use the normalized chain (cochain) complex
to compute the homology (cohomology) in Equation (\ref{Defi:limi}) (
Equation (\ref{Defi:liminvi})).
\end{Rmk}

There is a decreasing filtration of the chain complex
$(C_*(\C,F),d)$ given by

$$L^pC_n(\C,F)=\bigoplus_{\sigma \in {N\C}_n, deg(\sigma_n)\geq
p} F_{\sigma}.$$

It is straightforward that the triple $(C_*(\C,F),d,L^*)$ is a
filtered differential graded $\Z$-module, so it yields a spectral
sequence $(E^{*,*}_r,d_r)$ of cohomological type whose
differential $d_r$ has bidegree $(r,1-r)$. The $E^{*,*}_1$ page is
given by
$$
    E^{p,q}_1\simeq H^{p+q}(L^pC/L^{p-1}C)
\textit{.}$$ The differential graded $\Z$-module $L^pC/L^{p-1}C$
is in fact a simplicial abelian group because the face operators
$d_i$ and the degeneracy operators $s_i$ respect the filtration
$L^*$. The $n$-simplices are
$$
(L^pC/L^{p-1}C)_n=\bigoplus_{\sigma \in {N\C}_n, deg(\sigma_n)=p}
F_{\sigma}.
$$
Moreover, for each $p$, $L^pC/L^{p-1}C$ can be filtered again by
the condition $deg(\sigma_0)\leq p'$ to obtain a homological type
spectral sequence. Then arguing as above we obtain:
\begin{Prop}\label{spectral_sequence_1}
For a (decreasing) graded poset $\C$ and a functor $F:\C\rightarrow
\Ab$:
\begin{itemize}
\item There exists a cohomological type spectral sequence
$E^{*,*}_*$ with target $\limdir_* F$. \item There exists a
homological type spectral sequence $(E^p)_{*,*}^*$ with target the
column $E^{p,*}_1$ for each $p$.
\end{itemize}
\end{Prop}
Notice that the column $E^{p,*}_1$ is given by the cohomology of the
simplicial abelian group formed by the simplices that end on objects
of degree $p$, the column $(E^p)_{p',*}^1$ is given by the homology
of the simplicial abelian group formed by the simplices that end on
degree $p$ and begin on degree $p'$, and all the differentials in
the spectral sequences above are induced by the completely described
differential of $(C_*(\C,F),d)$.

As $\bigcup_p L^pC_n = C_n$ and $\bigcap_p L^pC_n = {0}$ for each
$n$ the spectral sequence $E^{*,*}_*$ converges weakly to its
target. In case the map $deg$ has a bounded image, i.e., when $\C$
is $N$ dimensional, the filtration $L^*$ is bounded below and
over, and so $E^{*,*}_*$ collapses after a finite number of pages.
The same assertions on weak converge and boundedness hold for the
spectral sequences $(E^p)_{*,*}^*$.

If we proceed in reverse order, i.e., filtrating first by the
degree of the beginning object and later by the degree of the
ending object, we obtain:

\begin{Prop}\label{spectral_sequence_2}
For a (decreasing) graded poset $\C$ and a functor $F:\C\rightarrow
\Ab$:
\begin{itemize}
\item There exists a homological type spectral sequence
$E_{*,*}^*$ with target $\limdir_* F$. \item There exists a
cohomological type spectral sequence $(E_p)^{*,*}_*$ with target the
column $E_{p,*}^1$ for each $p$.
\end{itemize}
\end{Prop}

If the degree function we take is increasing then the appropriate
conditions for the filtrations are $deg(\sigma_n)\leq p$ and
$deg(\sigma_0)\geq p'$, and the spectral sequences obtained in
Propositions \ref{spectral_sequence_1} and
\ref{spectral_sequence_2} are of homological (cohomological) type
instead of cohomological (homological) type.

\begin{table}
\begin{tabular}{|c|c|c|c|c|c|}
\hline
Complex & Degree    & First       & Second     & First ss. & Second ss. \\
        & function  & filtration  & filtration &           &            \\
\hline
$C_*(\C,F)$ & decreasing    & $deg(\sigma_n)\geq$ &  $deg(\sigma_0)\leq$ & cohomol. type &  homol. type \\
$C_*(\C,F)$ & decreasing    & $deg(\sigma_0)\leq$ &  $deg(\sigma_n)\geq$ & homol. type   &  cohomol. type \\
$C_*(\C,F)$ & increasing    & $deg(\sigma_n)\leq$ &  $deg(\sigma_0)\geq$ & homol. type   &  cohomol. type \\
$C_*(\C,F)$ & increasing    & $deg(\sigma_0)\geq$ &  $deg(\sigma_n)\leq$ & cohomol. type &  homol. type \\
$C^*(\C,F)$ & decreasing    & $deg(\sigma_n)\leq$ &  $deg(\sigma_0)\geq$ & homol. type   &  cohomol. type \\
$C^*(\C,F)$ & decreasing    & $deg(\sigma_0)\geq$ &  $deg(\sigma_n)\leq$ & cohomol. type &  homol. type \\
$C^*(\C,F)$ & increasing    & $deg(\sigma_n)\geq$ &  $deg(\sigma_0)\leq$ & cohomol. type &  homol. type \\
$C^*(\C,F)$ & increasing    & $deg(\sigma_0)\leq$ &  $deg(\sigma_n)\geq$ & homol. type   &  cohomol. type \\
\hline
\end{tabular}
\caption{Filtrations and spectral sequences
obtained\label{tabla_ss}}
\end{table}
For the case of the cochain complex $(C^*(\C,F),d)$ the choices for
the filtrations are $deg(\sigma_n)\leq p$ and $deg(\sigma_0)\geq p'$
for a decreasing degree function and $deg(\sigma_n)\geq p$ and
$deg(\sigma_0)\leq p'$ for an increasing one. Analogously we obtain
spectral sequences with target $\liminv^i F$ which columns in the
first page are computed by another spectral sequence.

Table \ref{tabla_ss} shows a summary of the types of the spectral
sequences for all the cases. The statements on weak convergence
and boundedness apply to any of the spectral sequences of the
table.

\begin{Rmk}\label{normalized_ss}
It is straightforward that normalizing (see Remark
\ref{simp_cosimp_normalization}) the simplicial (cosimplicial)
abelian groups that computes the page $1$ of the spectral sequences
above has the same effect as considering the spectral sequences of
the normalizations of $C_*(\C,F)$ ($C^*(\C,F)$).
\end{Rmk}


%% file: projective.tex
\section{Projective objects in
$\Ab^\P$.}\label{section_projective} Consider the abelian category
$\Ab^\P$ for some graded poset $\P$. In this section we shall
determine the projective objects in $\Ab^\P$. Recall that in $\Ab$
the projective objects are the free abelian groups. Along the rest
of the section $\P$ denotes a graded poset.

Suppose $F\in \Ab^\P$ is projective. How does $F$ look? Consider an
object $i_0$ of $\P$. We show that the quotient of $F(i_0)$ by the
images of the non-identity morphisms arriving to $i_0$ is free
abelian. To prove it, write

\begin{Defi}\label{defi_im}
$\im_F(i_0)=\sum_{i\stackrel{\alpha}\rightarrow i_0,\alpha\neq
1_{i_0}} \im F(\alpha)$ (or $\im_F(i_0)=0$ if the index set of the
sum is empty) and $\coker_F(i_0)=F(i_0)/\im_F(i_0)$.
\end{Defi}

It is straightforward that for a fixed object $i_0$ of $\P$ there is
a functor $$\coker_\cdot(i_0):\Ab^\P\rightarrow \Ab$$ which maps $F$
to $\coker_F(i_0)$. This functor is left adjoint to the skyscraper
functor $\Ab\rightarrow \Ab^\P$ which maps the abelian group $A_0$
to the functor $A:\P\rightarrow \Ab$ with values
\begin{numcases}{A(i)=}
A_0 & for $i=i_0$ \nonumber \\
0 & for $i\neq i_0$ \nonumber
\end{numcases}
on objects, and values
\begin{numcases}{A(\alpha)=}
1_{A_0} & for $\alpha=1_{i_0}$ \nonumber \\
0 & for $\alpha\neq 1_{i_0}$ \nonumber
\end{numcases}
on morphisms. As the skyscraper functor is exact we obtain by
\cite[Proposition 2.3.10]{weibel} that $\coker_{.}(i_0)$ preserves
projective objects:
\begin{Lem}\label{lem_proj_pseudo_graded_coker_projective}
Let $F:\P\rightarrow \Ab$ be a projective functor over a graded
poset $\P$. Then $\coker_F(i_0)$ is free abelian for every object
$i_0$ of $\P$.
\end{Lem}
This means that we can write
$$ F(i_0)=\im_F(i_0)\oplus\coker_F(i_0)$$
with $\coker_F(i_0)$ free abelian for every object $i_0$ of $\P$,
and also that
\begin{Ex}
For the category $\P$ with shape
$$\cdot \rightarrow \cdot $$
the functor $F:\P\rightarrow \Ab$ with values
$$\Z \stackrel{\times n}\rightarrow \Z $$ is not projective as $\coker_F$ on the right object equals
the non-free abelian group $\Z/n$.
\end{Ex}

Now that we know a little about the values that a projective
functor $F:\P\rightarrow \Ab$ takes on objects we can wonder about
the values $F(\alpha)$ for $\alpha\in \Hom(\P)$. Do they have any
special property? Recall that a feature of graded posets is that
there is at most one arrow between any two objects, and also that
\begin{Rmk}\label{im_on_graded}
If $\P$ is graded then for any object $i_0$ of $\P$
$$\im_F(i_0)=\sum_{i\stackrel{\alpha}\rightarrow i_0,deg(\alpha)=1} \im F(\alpha) $$
because every morphism factors as composition of morphisms of
degree $1$.
\end{Rmk}

We prove that the following property holds for $F$:

\begin{Defi}\label{property_pseudo-projective}
Let $F:\P\rightarrow \Ab$ be a functor over a graded poset $\P$ with
degree function $deg$. Given $d\geq 0$ we say that $F$ is
\emph{$d$-pseudo-projective} if for any object $i_0$ of $\P$ and $k$
different objects $i_j$ in $\P$, arrows $\alpha_j:i_j\rightarrow
i_0$ with $deg(\alpha_j)=d$, and $x_j\in F(i_j)$ $j=1,..,k$ such
that
$$\sum_{j=1,..,k}F(\alpha_j)(x_j)=0$$
we have that $x_j\in\im_F(i_j)$ $j=1,..,k$. If $F$ is
$d$-pseudo-projective for each $d\geq 0$ we call $F$
\emph{pseudo-projective}.
\end{Defi}

\begin{Rmk}\label{pseudo_injective_is injective}
In case $k=1$ and $\im_F(i_1)=0$ the condition states that
$F(\alpha_1)$ is a monomorphism. Notice that any functor is
$0$-pseudo-projective as the identity is a monomorphism.
\end{Rmk}

Before proving that projective functors $F$ over a graded poset
verify this property we define two functors $\coker_F$ and
$\coker_F'$ and natural transformations $\sigma$ and $\pi$ that fit
in the diagram
$$
\xymatrix{
             & F \ar@{=>}[d]^{\sigma} &\\
\coker_F' \ar@{=>}[r]^\pi & \coker_F \ar@{=>}[r] & 0 }
$$
for any functor $F:\P\rightarrow \Ab$ with $\P$ a graded poset. We
begin defining $\coker_F$. Because for every $\alpha:i_1\rightarrow
i_0$ holds that $F(\alpha)(\im_F(i_1))\leq \im_F(i_0)$ we can factor
$F(\alpha)$ as in the diagram
$$
\xymatrix{ F(i_1)\ar@{->>}[d]\ar[r]^{F(\alpha)} & F(i_0)\ar@{->>}[d]\\
\coker_F(i_1)\ar[r]^{\overline{F(\alpha)}} & \coker_F(i_0).
 }
$$
In fact, if $\alpha\neq 1_{i_1}$, then $\overline{F(\alpha)}\equiv
0$ by definition. Because the identity $1_{i_0}$ cannot be
factorized (by non-identity morphisms) in a graded poset then we
have a functor $\coker_F$ with value $\coker_F(i)$ on the object $i$
of $\P$ and which maps the non-identity morphisms to zero.
$\coker_F$ is a kind of ``discrete" functor. Also it is clear that
there exists a natural transformation $\sigma:F\Rightarrow \coker_F$
with $\sigma(i)$ the projection $F(i)\twoheadrightarrow
\coker_F(i)$.

Now we define $\coker_F'$ from $\coker_F$ in a similar way as free
diagrams are constructed. Let $\coker_F'$ be defined on objects by
$$\coker_F'(i_0)=\bigoplus_{\alpha:i\rightarrow i_0} \coker_F(i).$$
For $\beta\in \Hom(\P)$, $\beta:i_1\rightarrow i_0$,
$\coker_F'(\beta)$ is the only homomorphism which makes commute the
diagram
$$
\xymatrix{ \coker_F'(i_1)\ar[r]^{\coker_F'(\beta)} & \coker_F'(i_0) \\
\coker_F(i)\ar@{^{(}->}[u]\ar[r]^{1}& \coker_F(i)\ar@{^{(}->}[u]
 }
$$
for each $\alpha:i\rightarrow i_1$. In the bottom row of the
diagram, the direct summands $\coker_F(i)$ of $\coker_F'(i_1)$ and
$\coker_F'(i_0)$ correspond to $\alpha:i\rightarrow i_1$ and to the
composition $i\stackrel{\alpha}\rightarrow
i_1\stackrel{\beta}\rightarrow i_0$ respectively.

Then there exists a candidate to natural transformation
$\pi:\coker_F'\Rightarrow \coker_F$ which value $\pi(i)$ is the
projection $\pi(i):\coker_F'(i) \twoheadrightarrow \coker_F(i)$ onto
the direct summand corresponding to $1_i:i\rightarrow i$. Thus,
$\pi$ is a natural transformation if for every $\beta:i_1\rightarrow
i_0$ with $i_1\neq i_0$ the following diagram is commutative
$$
\xymatrix{\coker_F'(i_1)\ar@{->>}[d]^{\pi(i_1)}\ar[r]^{\coker_F'(\beta)} & \coker_F'(i_0)\ar@{->>}[d]^{\pi(i_0)}\\
\coker_F(i_1)\ar[r]^{0} & \coker_F(i_0).}
$$
It is clear that this square commutes if the identity $1_{i_0}$
cannot be factorized (by non-identity morphisms), and this holds
in a graded poset.

Now we have the commutative triangle
$$
\xymatrix{
             & F \ar@{=>}[d]^{\sigma} \ar@{==>}[dl]^\rho &\\
\coker_F' \ar@{=>}[r]^\pi & \coker_F \ar@{=>}[r] & 0 }
$$
where the natural transformation $\rho$ exists because $F$ is
projective. To prove that $F$ is $d$-pseudo-projective for some
$d\geq 0$ take an object $i_0$ of $\P$, $k$ objects $i_1,..,i_k$,
arrows $\alpha_j:i_j\rightarrow i_0$ with $deg(\alpha_j)=d$ and
elements $x_j\in F(i_j)$ for $j=1,..,k$ such that
$$\sum_{j=1,..,k}F(\alpha_j)(x_j)=0.$$

To visualize what is going on consider the diagram above near
$i_0$ for $k=2$
\begin{footnotesize}
$$
\xymatrix{
&&&  F(i_1)\ar@{-->}[ddlll]\ar@{->>}[dd]\ar[rd]^{F(\alpha_1)}& & F(i_2)\ar@<-6pt>@{-->}[ddlll]\ar@{->>}[dd]\ar[dl]^{F(\alpha_2)}\\
&&&  & F(i_0)\ar@<4pt>@{-->}[ddlll]\ar@{->>}[dd] &\\
\coker_F'(i_1)\ar[rd]      ^{\coker_F'(\alpha_1)} & &
\coker_F'(i_2)\ar[ld]_{\coker_F'(\alpha_2)}&
 \coker_F(i_1)\ar[rd]^{0}& & \coker_F(i_2)\ar[dl]_{0}\\
& \coker_F'(i_0)\ar@{->>}[rrr]^{\pi(i_0)} &&
 & \coker_F(i_0) &
}
$$
\end{footnotesize}
where $\pi$ is not drawn completely for clarity. Recall that we
are supposing that $\{x_1,..,x_k\}$ is such that
$\sum_{j=1,..,k}F(\alpha_j)(x_j)=0$. Then
$$
0=\rho(i_0)(0)=\sum_{j=1,..,k}\rho(i_0)(F(\alpha_j)(x_j))
 = \sum_{j=1,..,k}\coker_F'(\alpha_j)(\rho(i_j)(x_j)).
$$

Now consider the projection $p_{j_0}$ for $j_0\in\{1,..,k\}$ from
$\coker_F'(i_0)$ onto the direct summand
$\coker_F(i_{j_0})\hookrightarrow \coker_F'(i_0)$ which corresponds
to $\alpha_{j_0}:i_{j_0}\rightarrow i_0$
$$\coker_F'(i_0)\stackrel{p_{j_0}}\twoheadrightarrow \coker_F(i_{j_0}).$$
Then
\begin{equation}\label{projiszero}
 0=p_{j_0}(0)=p_{j_0}(\rho(i_0)(0))=\sum_{j=1,..,k}
p_{j_0}(\coker_F'(\alpha_j)(\rho(i_j)(x_j))).
\end{equation}

For any $y=\bigoplus_{\alpha:i\rightarrow i_j}
y_\alpha\in\coker_F'(y_j)$
$$p_{j_0}(\coker_F'(\alpha_j)(y))=\sum_{\alpha:i\rightarrow i_j,\alpha_j\circ\alpha=\alpha_{j_0}}
y_\alpha.$$

So if $y_j=\rho(i_j)(x_j)=\bigoplus_{\alpha:i\rightarrow i_j}
y_{j,\alpha} \in \coker_F'(i_j)$ then
$$p_{j_0}(\coker_F'(\alpha_j)(\rho(i_j)(x_j)))=\sum_{\alpha:i\rightarrow i_j,\alpha_j\circ\alpha=\alpha_{j_0}}
y_{j,\alpha}.$$

This last sum runs over $\alpha:i_{j_0}\rightarrow i_j$ such that
the following triangle commutes
$$
\xymatrix{ i_{j_0}\ar[rd]^{\alpha}\ar[rr]^{\alpha_{j_0}} && i_0.\\
& i_j \ar[ru]^{\alpha_j} &}
$$
Because we are in a graded poset and $deg(i_j)=d$ for each
$j=1,..,k$ then the only chance is $i_j=i_{j_0}$ and
$\alpha=1_{i_{j_0}}$. Because the objects $i_1,..,i_k$ are
different this implies that $j=j_0$ too. Thus
\begin{numcases}{p_{j_0}(\coker_F'(\alpha_j)(\rho(i_j)(x_j)))=}
y_{j_0,1_{i_{j_0}}} & for $j=j_0$ \nonumber \\
0 & for $j\neq j_0$ \nonumber
\end{numcases}
and Equation (\ref{projiszero}) becomes
$$
0=p_{j_0}(0)=y_{j_0,1_{i_{j_0}}}.
$$
Notice now that $y_{j_0,1_{i_{j_0}}}$ is the evaluation of
$\pi(i_{j_0})$ on $y_{j_0}=\rho(i_{j_0})(x_{j_0})$ and then
$$
0=y_{j_0,1_{i_{j_0}}}=\pi(i_{j_0})(\rho(i_{j_0})(x_{j_0}))=\sigma_{i_{j_0}}(x_{j_0}).
$$
This last equation means that $x_{j_0}$ goes to zero by the
projection $F(i_{j_0})\twoheadrightarrow
\coker_F(i_{j_0})=F(i_{j_0})/\im_F(i_{j_0})$, and then
$$x_{i_{j_0}}\in \im_F(i_{j_0}).$$
As $j_0$ was arbitrary this completes the proof of

\begin{Lem}\label{lem_proj_pseudo_graded}
Let $F:\P\rightarrow \Ab$ be a projective functor over a graded
poset $\P$. Then $F$ is pseudo-projective.
\end{Lem}

\begin{Ex}
For the category $\P$ with shape
$$\cdot \rightarrow \cdot $$
the functor $F:\P\rightarrow \Ab$ with values
$$\Z \stackrel{red_n}\rightarrow \Z/n $$ is not projective as
$red_n$ is not injective, in spite of the $\coker_F$'s are $\Z$ and
$0$, which are free abelian.
\end{Ex}
Till now we have obtained (Lemmas
\ref{lem_proj_pseudo_graded_coker_projective} and
\ref{lem_proj_pseudo_graded}) that projective functors
$\P\rightarrow \Ab$ are pseudo-projective and have $\coker_F(i_0)$
projective for any object $i_0$. In fact, as the next theorem shows,
the restriction we did to graded posets is worthwhile:

\begin{Thm}\label{pro_projective_bounded_graded}
Let $F:\P\rightarrow \Ab$ be a functor over a bounded below graded
poset $\P$. Then $F$ is projective if and only if
\begin{enumerate}
\item for any object $i_0$ of $\P$ $\coker_F(i_0)$ is a free abelian group.
\label{Defi_preprojective_1}\item $F$ is pseudo-projective.
\label{Defi_preprojective_2}
\end{enumerate}
\end{Thm}
\begin{proof}
It remains to prove that a functor $F$ satisfying the conditions in
the statement is projective. We can assume that the degree function
$deg$ on $\P$ is increasing and takes values $\{0,1,2,3,...\}$, and
that $\Ob_0(\P)\neq \emptyset$.

To see that $F$ is proyective in $\Ab^\P$, given a diagram of
functors with exact row as shown, we must find a natural
transformation $\rho:F\Rightarrow A$ making the diagram
commutative:
$$
\xymatrix{
             & F \ar@{=>}[d]^\sigma \ar@{==>}[dl]^\rho &\\
A \ar@{=>}[r]^\pi & B \ar@{=>}[r]                           & 0. }
$$
We define $\rho$ inductively, beginning on objects of degree $0$
and successively on object of degrees $1,2,3,..$.

So take $i_0\in \Ob_0(\P)$ of degree $0$, and restrict to the
diagram in $\Ab$ over $i_0$. By hypothesis
(\ref{Defi_preprojective_1}) in the statement,
as $\im_F(i_0)=0$, $F(i_0)=\coker_F(i_0)$ is free abelian. So we can
close the following triangle with a homomorphism $\rho(i_0)$
$$
\xymatrix{
             & F(i_0) \ar[d]^{\sigma(i_0)} \ar@{-->}[dl]_{\rho(i_0)} &\\
A(i_0) \ar[r]^{\pi(i_0)} & B(i_0) \ar[r] & 0. }
$$
As there are no arrows between degree $0$ objects we do not worry
about $\rho$ being a natural transformation. Now suppose that we
have defined $\rho$ on all objects of $\P$ of degree less than $n$
($n\geq 1$), and that the restriction of $\rho$ to the full
subcategory generated by these objects is a natural transformation
and verifies $\pi\circ \rho=\sigma$.

The next step is to define $\rho$ on degree $n$ objects. So take
$i_0\in \Ob_n(\P)$ and consider the splitting
$$F(i_0)=\im_F(i_0)\oplus\coker_F(i_0)$$
where
$$\im_F(i_0)=\sum_{i\stackrel{\alpha}\rightarrow i_0,deg(\alpha)=1} \im F(\alpha).$$
To define $\rho(i_0)$ such that it makes commutative the diagram
$$
\xymatrix{
             & \im_F(i_0)\oplus\coker_F(i_0) \ar[d]^{\sigma(i_0)} \ar@{-->}[dl]_{\rho(i_0)} &\\
A(i_0) \ar[r]^{\pi(i_0)} & B(i_0) \ar[r] & 0,}
$$
we define it on $\im_F(i_0)$ and $\coker_F(i_0)$ separately. For
$\coker_F(i_0)$, as it is a free abelian group, we define it by any
homomorphism that makes commutative the diagram above when
restricted to $\coker_F(i_0)$. For $\im_F(i_0)$ take
$x=\sum_{j=1,..,k} F(\alpha_j)(x_j)$ where $\{i_1,..,i_k\}$ are $k$
different objects, $\alpha_j:i_j\rightarrow i_0$, $deg(\alpha_j)=1$
and $x_j\in F(i_j)$ for $j=1,..,k$ (see Remark \ref{im_on_graded}).
Then define
$$
\rho(i_0)(x)=\sum_{j=1,..,k} (A(\alpha_j)\circ\rho(i_j))(x_j).
$$
To check that $\rho(i_0)(x)$ does not depend on the choice of the
$i_j$'s, $\alpha_j$'s and $x_j$'s we have to prove that
$$
\sum_{j=1,..,k} F(\alpha_j)(x_j)=0\Rightarrow \sum_{j=1,..,k}
(A(\alpha_j)\circ \rho(i_j))(x_j)=0.
$$

So suppose that
\begin{equation}\label{equ_pre_pro_0}
\sum_{j=1,..,k}F(\alpha_j)(x_j)=0.
\end{equation}
Then using that $F$ is $1$-pseudo-projective and Remark
\ref{im_on_graded} we obtain objects $i_{j,j'}$, arrows
$\alpha_{j,j'}$ of degree $1$, and elements $x_{j,j'}$ for
$j=1,..,k$, $j'=1,..,k_j$ such that
\begin{equation}\label{equ_pre pro_1}
\sum_{j'=1,..,k_j}F(\alpha_{j,j'})(x_{j,j'})=x_j
\end{equation}
for every $j\in\{1,..,k\}$. Notice that possibly not all the
objects $i_{j,j'}$ are different.  Replacing Equation
(\ref{equ_pre pro_1}) in Equation (\ref{equ_pre_pro_0}) we obtain
\begin{equation}\label{equ_pre_pro_2}
\sum_{\textit{$j=1,..,k$, $j'=1,..,k_j$}}
F(\alpha_j\circ\alpha_{j,j'})(x_{j,j'})=0.
\end{equation}
Because in a graded poset there is at most one arrow between two
objects, the condition $i_{j,j'}=i_{j',{j'}'}=i$ implies
$\alpha_j\circ\alpha_{j,j'}=\alpha_{j'}\circ\alpha_{j',{j'}'}:i\rightarrow
i_0$. So, considering objects $i$ in $\P$, we can rewrite
(\ref{equ_pre_pro_2}) as
\begin{equation}\label{equ pre pro 3}
\sum_{i\in \Ob(\P)} F(\alpha_j\circ\alpha_{j,j'})(\sum_{j,j' |
i_{j,j'}=i} x_{j,j'})=0.
\end{equation}

Call $\{i'_1,..,i'_m\}=\{i_{j,j'}|\textit{$j=1,..,k$,
$j'=1,..,k_j$}\}$ where these sets have $m$ elements. Call
$\beta_l=\alpha_j\circ \alpha_{j,j'}$ if $i'_l=i_{j,j'}$ and
$y_l=\sum_{j,j' | i_{j,j'}=i'_l} x_{j,j'}$ for $l=1,..,m$. Notice
that $deg(\beta_l)=2$ for each $l$. Then Equation (\ref{equ pre
pro 3}) becomes
\begin{equation}\label{equ_pre pro 4}
\sum_{l=1,..,m} F(\beta_l)(y_l)=0.
\end{equation}

Now we repeat the same argument: applying that $F$ is
$2$-pseudo-projective and the Remark \ref{im_on_graded} to
Equation (\ref{equ_pre pro 4}) we obtain objects $i'_{l,l'}$,
arrows $\beta_{l,l'}$ of degree $1$, and elements $y_{l,l'}$ for
$l=1,..,m$, $l'=1,..,k'_l$ such that
\begin{equation}\label{equ_pre pro_5}
\sum_{l'=1,..,k'_l}F(\beta_{l,l'})(y_{l,l'})=y_l
\end{equation}
for every $l\in\{1,..,m\}$. Substituting (\ref{equ_pre pro_5}) in
(\ref{equ_pre pro 4})
$$
\sum_{\textit{$l=1,..,m$,
$l'=1,..,k'_l$}}F(\beta_l\circ\beta_{l,l'})(y_{l,l'})=0.
$$
Now proceed as before regrouping the terms in this last equation.

In a finite number of steps, after a regrouping of terms as above,
we find objects $i''_s$, arrows $\gamma_s$, and elements $z_s$ of
degree $0$ for $s=1,..,r$ which verify an equation
\begin{equation}\label{equ_pre pro_6}
\sum_{s=1,..,r}F(\gamma_s)(z_s)=0.
\end{equation}
Then pseudo-injectivity gives that $z_s\in \im_F(i''_s)$ for each
$s$. As $deg(i''_s)=0$ then $\im_F(i''_s)=0$ and so $z_s=0$ (notice
that $z_s=0$ for $s=1,..,r$ does not imply $x_j=0$ for any $j$).

Recall that we want to prove that
\begin{equation}\label{equ pre pro 7}
\sum_{j=1,..,k} (A(\alpha_j)\circ\rho(i_j))(x_j)=0.
\end{equation}
Substituting (\ref{equ_pre pro_1}) in
$\sum_{j=1,..,k}(A(\alpha_j)\circ\rho(i_j))(x_j)$ we obtain
\begin{align}
\sum_{j=1,..,k}(A(\alpha_j)\circ\rho(i_j))(x_j))&=\sum_{j=1,..,k}\sum_{j'=1,..,k_j}(A(\alpha_j)\circ\rho(i_j)\circ F(\alpha_{j,j'}))(x_{j,j'})\notag\\
&=\sum_{j=1,..,k}\sum_{j'=1,..,k_j}(A(\alpha_j)\circ
A(\alpha_{j,j'})\circ\rho(i_{j,j'}) )(x_{j,j'})\notag\\
&=\sum_{j=1,..,k}\sum_{j'=1,..,k_j}(A(\alpha_j\circ
\alpha_{j,j'})\circ\rho(i_{j,j'}) )(x_{j,j'}),\notag
\end{align}
as $\rho$ is natural up to degree less than $n$. Then regrouping
terms
\begin{align}
\sum_{j=1,..,k}\sum_{j'=1,..,k_j}(A(\alpha_j\circ
\alpha_{j,j'})\circ\rho(i_{j,j'}) )(x_{j,j'})&=\sum_{i\in
\Ob(\P)}(A(\alpha_j\circ \alpha_{j,j'})\circ
\rho(i_{j,j'}))(\sum_{j,j' | i_{j,j'}=i} x_{j,j'})\notag\\
&=\sum_{l=1,..,m}(A(\beta_l)\circ \rho(i'_l))(y_l).\notag
\end{align}

Then, after a finite number of steps, we obtain
$$
\sum_{j=1,..,k}(A(\alpha_j)\circ\rho(i_j))(x_j))=\sum_{\textit{$s=1,..,r$}}(A(\gamma_s)\circ\rho(i''_s))(z_s)=0
$$
as $z_s=0$ for each $z=1,..,r$.

So we have checked that $\rho(i_0)(x)$ does not depend on the choice
of $i_j$, $\alpha_j$ and $x_j$. It is straightforward that
$\rho(i_0)$ on $\im_F(i_0)$ defined in this way is a homomorphism of
abelian groups.

It remains to prove that $\pi(i_0)\rho(i_0)=\sigma(i_0)$ when
restricted to $\im_F(i_0)$. So take
$x=\sum_{j=1,..,k}F(\alpha_j)(x_j)$ in $\im_F(i_0)$. Then
\begin{align}
\pi(i_0)(\rho(i_0)(x)) =&\sum_{j=1,..,k}(\pi(i_0)\circ A(\alpha_j)\circ \rho(i_j))(x_j)  \notag \\
=&\sum_{j=1,..,k}(B(\alpha_j)\circ\pi(i_j)\circ\rho(i_j))(x_j) \text{, $\pi$ is a natural transformation}\notag\\
=&\sum_{j=1,..,k}(B(\alpha_j)\circ\sigma(i_j))(x_j) \text{, by the inductive hypothesis} \notag\\
=&\sum_{j=1,..,k}(\sigma(i_0)\circ F(\alpha_j))(x_j) \text{, $\sigma$ is a natural transformation} \notag\\
=&\sigma(i_0)(x) \notag
\end{align}

Defining $\rho(i_0)$ in this way for every $i_0\in \Ob_n(\P)$ we
have now $\rho$ defined on all objects of $\P$ of degree less or
equal than $n$. Finally, to complete the inductive step we have to
prove that $\rho$ restricted to the full subcategory over these
objects is a natural transformation. Take $\alpha:i\rightarrow
i_0$ in this full subcategory. If the degree of $i_0$ is less than
$n$ then the commutativity of
$$
\xymatrix{
F(i)\ar[r]^{F(\alpha)}\ar[d]^{\rho(i)}   &   F(i_0)\ar[d]^{\rho(i_0)} \\
A(i)\ar[r]^{A(\alpha)}                  &   A(i_0) \\
}
$$
is granted by the inductive hypothesis. Suppose that the degree of
$i_0$ is $n$. Take $x'\in F(i)$. Because $\P$ is graded there
exists $\alpha_1:i_1\rightarrow i_0$ of degree $1$ and
$\alpha':i\rightarrow i_1$ such that
$\alpha=\alpha_1\circ\alpha'$:
$$
\xymatrix{ i\ar[rd]^{\alpha'}\ar[rr]^{\alpha} && i_0.\\
& i_1 \ar[ru]^{\alpha_1} &}
$$
Write $x=F(\alpha)(x')=F(\alpha_1)(x_1)$ where $x_1=F(\alpha')(x')$.
Then, by definition of $\rho(i_0)$ on $\im_F(i_0)$,
\begin{align}
\rho(i_0)(x)&=(A(\alpha_1)\circ\rho(i_1))(x_1)\notag\\
&=(A(\alpha_1)\circ\rho(i_1))(F(\alpha')(x'))\notag\\
&=(A(\alpha_1)\circ\rho(i_1)\circ F(\alpha'))(x')\notag\\
&=(A(\alpha_1 \circ \alpha')\circ \rho(i))(x')\text{, $\rho$ is natural up to degree less than $n$ }\notag\\
&=(A(\alpha)\circ \rho(i))(x')\notag
\end{align}
and so the diagram commutes.
\end{proof}

\begin{Rmk}\label{remark:neccesary_bounded}
As the following example shows the condition of lower \linebreak
boundedness of $\P$ in Theorem \ref{pro_projective_bounded_graded}
cannot be dropped:

Consider the inverse `telescope category' $\P$ with shape
$$
...\rightarrow \cdot \rightarrow \cdot \rightarrow \cdot
$$
It is a graded poset which is not bounded below. Consider the
functor of constant value $\Z/p$, $c_{\Z/p}:\P\rightarrow \Ab$:
$$
...\rightarrow \Z/p \rightarrow \Z/p \rightarrow \Z/p
$$
It is straightforward that it satisfies the conditions in the
theorem as all the cokernels are zero and all the arrows are
injective. But it is not a projective object of $\Ab^\P$ because, in
that case, the adjoint pair $\limdir:\Ab^\P\leftrightarrow
\Ab:\Delta$ would give that $\Z/p$ is projective in $\Ab$ ( see
\cite[3.2, Ex7]{cohn} or \cite[Proposition 2.3.10]{weibel}).
\end{Rmk}

%

This theorem yields the following examples. The degree functions
$deg$ for the bounded below graded posets appearing in the examples
are indicated by subscripts $i_{deg(i)}$ on the objects $i$ of $\P$
and take values $\{0,1,2,3,...\}$.

\begin{Ex}\label{examples_projective}
For the `pushout category' $\P$ with shape
$$\xymatrix {a_0 \ar[r]^{f}\ar[d]^{g} & b_1 \\ c_1}$$
a functor $F:\P\rightarrow \Ab$ is projective if and only if
\begin{itemize}
\item $F(a)$, $F(b)/{\im F(f)}$ and $F(c)/{\im F(g)}$ are free
abelian. \item $F(f)$ and $F(g)$ are monomorphisms.
\end{itemize}

%
For the `telescope category' $\P$ with shape
$$ \xymatrix{a_0 \ar[r]^{f_1} & a_1 \ar[r]^{f_2} & a_2 \ar[r]^{f_3} &  a_3 \ar[r]^{f_4} & a_4 ...}$$
a functor $F:\P\rightarrow \Ab$ is projective if and only if
\begin{itemize}
\item $F(a_0)$ is free abelian. \item $F(a_i)/\im F(f_i)$ is free
abelian, $F(f_{i}\circ f_{i-1} \circ ..\circ f_0)$ is a
monomorphism and $\ker F(f_{i}\circ f_{i-1} \circ ..\circ
f_{i-d})\subseteq \im F(f_{i-d-1})$ for $d=0,1,..,i-1$  for each
$i=1,2,3,4,...$.
\end{itemize}
\end{Ex}

%% file: limacyI.tex
\section{Pseudo-projectivity}\label{section_Pseudo-projectivity}
Consider a functor $F:\P\rightarrow \Ab$ over a graded poset $\P$.
In this section we find conditions on $F$ such that $\limdir_i F=0$
for $i\geq 1$. We fix the following notation
\begin{Defi}\label{Defi_direct_acyclic}
Let $\P$ be a graded poset and $F:\P\rightarrow \Ab$. We say that
$F$ is \emph{$\limdir$-acyclic} if $\limdir_i F=0$ for $i\geq 1$.
\end{Defi}

Recall that for projective objects it holds that any left derived
functor vanishes.
%
Because being $\limdir$-acyclic is clearly weaker that being
projective we can wonder if is it possible to weaken the hypothesis
of projectiveness keeping the thesis of $\limdir$-acyclicity. It
turns out that pseudo-projectiveness gives an appropriate weaker
condition:

\begin{Thm}\label{pro_acyclic_graded}
Let $F:\P\rightarrow \Ab$ be a pseudo-projective functor over a
bounded below graded poset $\P$. Then $F$ is $\limdir$-acyclic.
\end{Thm}
\begin{proof}
We can suppose that the degree function $deg$ on $\P$ is
increasing and takes values $\{0,1,2,3,...\}$, and that
$\Ob_0(\P)\neq \emptyset$. To compute $\limdir_t F$ we use the
(normalized, Remark \ref{normalized_ss}) spectral sequences
corresponding to the third row of Table \ref{tabla_ss} in Chapter
\ref{spectral}. That is, we first filter by the degree of the end
object of each simplex to obtain a homological type spectral
sequence $E^*_{*,*}$. To compute the column $E^1_{p,*}$ we filter
by the degree of the initial object of each object to obtain
cohomological type spectral sequences $(E_p)_*^{*,*}$.

Fix $t\geq 1$. Notice that to prove that $\limdir_t F=0$ is enough
to show that $E^1_{p,t-p}$ is zero for every $p$. The
contributions to $E^1_{p,t-p}$ come from
$(E_p)_{\infty}^{p',p-p'-t}$ for $p'\leq p-t$ (we are using
normalized (Remark \ref{normalized_ss}) spectral sequences). We
prove that
$$
(E_p)_{r}^{p',p-p'-t}=0
$$
if $r$ is big enough for each $p$ and $p'\leq p-t$. This implies
that $\limdir_t F=0$.

Consider the increasing filtration $L^*$ of $C_*(\P,F)$ that gives
rise to the spectral sequence $E^*_{*,*}$. The $n$-simplices are
$$
L^p_n=L^pC_n(\P,F)=\bigoplus_{\sigma \in {N\P}_n,
deg(\sigma_n)\leq p} F_{\sigma}.
$$
For each $p$ we have a decreasing filtration $M_p^*$ of the
quotient $L^p/L^{p-1}$ that gives rise to the spectral sequence
$(E_p)_*^{*,*}$ and which $n$-simplices are
$$
(M_p)^{p'}_n=\bigoplus_{\sigma \in {N\P}_n, deg(\sigma_0)\geq p',
deg(\sigma_n)=p} F_{\sigma}.
$$
For $p'\leq p-t$ the abelian group $(E_p)_r^{p',q'}$ at the
$t=-(p'+q')+p$ simplices is given by
$$
(E_p)_r^{p',q'}= (M_p)^{p'}_t\cap
d^{-1}((M_p)^{p'+r}_{t-1})/(M_p)^{p'+1}_t\cap
d^{-1}((M_p)^{p'+r}_{t-1})+(M_p)^{p'}_t\cap
d((M_p)^{p'-r+1}_{t+1})
$$
where $d$ is the differential of the quotient $L^p/L^{p-1}$
restricted to the subgroups of the filtration $(M_p)^*$. For
$r>p-p'-(t-1)$ there are not $(t-1)$-simplices beginning in degree
at least $p'+r>p-(t-1)$ and ending in degree $p$, i.e.,
$(M_p)^{p'+r}_{t-1}=0$. Because $\P$ is bounded below  for $r$ big
enough $(M_p)^{p'-r+1}_{t+1}=(M_p)^0_{t+1}=(L^p/L^{p-1})_{t+1}$,
i.e., $(M_p)^{p'-r+1}_{t+1}$ equals all the $(t+1)$-simplices that
end on degree $p$. Thus there exists $r$ such that
\begin{equation}
\label{E_p_r_big_dir} (E_p)_r^{p',q'}= (M_p)^{p'}_t\cap
d^{-1}(0)/(M_p)^{p'+1}_t\cap d^{-1}(0)+(M_p)^{p'}_t\cap
d((M_p)^0_{t+1}).
\end{equation}
Fix such an $r$ and take $[x]\in(E_p)_r^{p',q'}$ where
\begin{equation}
\label{equ_exp_x_dir} x=\bigoplus_{\sigma \in {N\P}_t,
deg(\sigma_0)\geq p' , deg(\sigma_t)=p} x_{\sigma}
\end{equation}
and $d(x)=0$. Notice that by definition there is just a finite
number of summands $x_\sigma\neq 0$ in the expression
(\ref{equ_exp_x_dir}) for $x$. We prove that $[x]=0$ in three
steps:

\textbf{Step 1:} In this first step we find a representative $x'$
for $[x]$
$$
x'=\bigoplus_{\sigma \in {N\P}_t, deg(\sigma_0)\geq p' ,
deg(\sigma_t)=p} x'_{\sigma} $$  such that $\deg(\alpha_1)=1$ for
every $\sigma=\xymatrix{ \sigma_0 \ar[r]^{\alpha_1} & \sigma_1
\ar[r]^{\alpha_2}&...\ar[r]^{\alpha_{t-1}}&
\sigma_{t-1}\ar[r]^{\alpha_t} & \sigma_t}$ with $x'_\sigma\neq 0$.

Take $\sigma$ such that $x_\sigma\neq 0$ and suppose that
$deg(\alpha_1)>1$, i.e., $deg(\sigma_0)<deg(\sigma_1)-1$. Then, as
in a graded poset every morphism factors as composition of degree
$1$ morphisms, there exists an object $\sigma_*$ of degree
$deg(\sigma_0)<deg(\sigma_*)<deg(\sigma_1)$ and arrows
$\beta_1:\sigma_0\rightarrow \sigma_*$ and
$\beta_2:\sigma_*\rightarrow \sigma_1$ with $\alpha_1=\beta_2\circ
\beta_1$.
$$\xymatrix{& \sigma_* \ar[rd]^{\beta_2} & \\
\sigma_0 \ar[rr]^{\alpha_1} \ar[ru]^{\beta_1}& & \sigma_1. }$$

Call $\tilde{\sigma}$ to the $(t+1)$-simplex $\sigma=\xymatrix{
\sigma_0 \ar[r]^{\beta_1} & \sigma_* \ar[r]^{\beta_1} & \sigma_1
\ar[r]^{\alpha_2}&...\ar[r]^{\alpha_{t-1}}&
\sigma_{t-1}\ar[r]^{\alpha_t} & \sigma_t}$ and consider the
$(t+1)$-chain of $(M_p)^0_{t+1}$
$y=i_{\tilde{\sigma}}(-x_\sigma)$. Its differential in
$L^p/L^{p-1}$ equals
$$d(y)=d_0(y)-d_1(y)+\sum_{i=2,..,t}(-1)^i d_i(y)=d_0(y)+i_\sigma(x_\sigma)+\sum_{i=2,..,t}(-1)^i d_i(y).$$
Notice that the first morphisms appearing in the simplices
$d_0(\tilde{\sigma})$ and $d_i(\tilde{\sigma})$ for $i=2,..,t$
have degree $deg(\beta_2)$ and $deg(\beta_1)$ respectively, which
are strictly less than $deg(\alpha_1)$. Also notice that $d(y)\in
(M_p)^{p'}_t\cap d((M_p)^0_{t+1})$ (which is  zero in Equation
(\ref{E_p_r_big_dir}) ).

Taking the (finite) sum of the chains $y$ for each term $x_\sigma$
we find that $[x]=[x']$ where
$$
x'=\bigoplus_{\sigma \in {N\P}_t, deg(\sigma_0)\geq p' ,
deg(\sigma_t)=p} x'_{\sigma} $$ and the maximum of the degrees of
the morphisms $\alpha_1$ of the simplices
$$\sigma=\xymatrix{ \sigma_0 \ar[r]^{\alpha_1} & \sigma_1
\ar[r]^{\alpha_2}&...\ar[r]^{\alpha_{t-1}}&
\sigma_{t-1}\ar[r]^{\alpha_t} & \sigma_t}$$ with $x'_\sigma \neq
0$ is smaller than this maximum computed for $x$. So repeating
this process a finite number of times we find a representative as
wished. For simplicity we write also $x$ for this representative.

\textbf{Step 2:} By Step $1$ we can suppose that
$\deg(\alpha_1)=1$ for every $$\sigma=\xymatrix{ \sigma_0
\ar[r]^{\alpha_1} & \sigma_1
\ar[r]^{\alpha_2}&...\ar[r]^{\alpha_{t-1}}&
\sigma_{t-1}\ar[r]^{\alpha_t} & \sigma_t}$$ with $x_\sigma\neq 0$.
Now our objective is to find a representative $x'$ for $[x]$
$$
x'=\bigoplus_{\sigma \in {N\P}_t, deg(\sigma_0)=p' ,
deg(\sigma_t)=p} x'_{\sigma},$$ i.e., such that the expression for
$x'$ runs over simplices $\sigma$ with begin in degree $p'$. Begin
writing $x$ as
$$
x=\bigoplus_{i=p',..,p-t} x_i
$$
where
$$
x_i=\bigoplus_{\sigma \in {N\P}_t, deg(\sigma_0)=i,
deg(\sigma_t)=p} x_\sigma.
$$
Notice that the index $i$ just goes to $p-t$ (and not to $p$)
because we are using normalized (Remark \ref{normalized_ss})
spectral sequences. Now we prove
\begin{Claim}\label{claim_step_2}
For each $i$ from $i=p-t$ to $i=p'$ there exists a representative
$x'_i$ for $[x]$
$$
x'_i=\bigoplus_{\sigma \in {N\P}_t,i\geq deg(\sigma_0)\geq p' ,
deg(\sigma_t)=p} (x'_i)_{\sigma} $$ such that
\begin{equation}
\label{condition_claim_step_2} \textit{$(x'_i)_\sigma\neq 0$ and
$deg(\sigma_0)<i$ imply $deg(\alpha_1)=1$.}
\end{equation}
\end{Claim}
Notice that taking $i=p'$ in the claim, the step $2$ is finished.
The case $i=p-t$ in the claim is fulfilled taking $x'_{p-t}=x$ (by
step $1$). Suppose the statement of the claim holds for $i$. Then
we prove it for $i-1$. We have $x'_i$ such that
$$ x'_i=\bigoplus_{\sigma \in {N\P}_t,i\geq deg(\sigma_0)\geq p' ,
deg(\sigma_t)=p} (x'_i)_{\sigma}, $$ $d(x'_i)=0$ and $[x]=[x'_i]$.
The differential $d$ on $L^p/L^{p-1}$ restricts to
$$d:(M_p)^{p'}_t\rightarrow (M_p)^{p'}_{t-1}$$
and carries $z\in F_\sigma\hookrightarrow \bigoplus_{\sigma \in
{N\P}_t, deg(\sigma_0)\geq p', deg(\sigma_t)=p}
F_{\sigma}=(M_p)^{p'}_t$ to
$$
d(z)=\sum_{j=0,1,..,t-1} (-1)^j d_j(z)
$$
with $d_j(z)\in F_{d_j(\sigma)}\hookrightarrow (M_p)^{p'}_{t-1}$.
Notice that the initial object of $d_j(\sigma)$ is $\sigma_1$ for
$j=0$ and $\sigma_0$ for $j=1,..,t-1$. Also notice that the final
object of $d_j(\sigma)$ is $\sigma_t$ for $j=0,..,t-1$.

By hypothesis $d(x'_i)=0$. So for every $\epsilon\in {N\P}_{t-1}$
with $deg(\epsilon_0)\geq p'$ and $deg(\epsilon_{t-1})=p$ we can
apply the projection
$$\pi_\epsilon: (M_p)^{p'}_{t-1} \twoheadrightarrow F_\epsilon$$ and
obtain $\pi_\epsilon(d(x'_i))=0$. If $deg(\epsilon_0)>i$ then the
remarks on the differential above and condition
(\ref{condition_claim_step_2}) imply that
$$\pi_\epsilon(d(x'_i))=\sum_{\sigma \in{N\P}_t, deg(\sigma_0)=i, d_0(\sigma)=\epsilon}
F(\alpha_1)((x'_i)_{\sigma})
$$
and thus
\begin{equation}\label{equ_epsilon}
0=\sum_{\sigma \in{N\P}_t, deg(\sigma_0)=i, d_0(\sigma)=\epsilon}
F(\alpha_1)((x'_i)_{\sigma})
\end{equation}
for each $\epsilon\in {N\P}_{t-1}$ with $deg(\epsilon_0)>i$ and
$deg(\epsilon_{t-1})=p$. Notice that each summand $(x'_i)_\sigma$
with $\sigma \in {N\P}_t$, $deg(\sigma_0)=i$ and $deg(\sigma)=p$
appears in one and just one equation as (\ref{equ_epsilon}) (take
$\epsilon=d_0(\sigma)$).

Fix an $\epsilon\in {N\P}_{t-1}$ with $deg(\epsilon_0)>i $ and
$deg(\epsilon_{t-1})=p$ and consider the associated Equation
(\ref{equ_epsilon}). Then, as $F$ is
$(i-deg(\epsilon_0))$-pseudo-projective, $(x'_i)_\sigma\in
\im_F(\sigma_0)$ for every $\sigma \in{N\P}_t$ with
$deg(\sigma_0)=i$ and $d_0(\sigma)=\epsilon$. This means that for
every such a $\sigma$ there exists $k_\sigma$ objects of degree
$(i-1)$, namely $i_\sigma^1,..,i_\sigma^{k_\sigma}$, arrows
$\beta_\sigma^j:i_\sigma^j\rightarrow \sigma_0$ and elements
$x_\sigma^j\in F(i_\sigma^j)$ for $j=1,..,k_\sigma$ such that
\begin{equation}\label{xenim}
(x'_i)_\sigma=\sum_{j=1,..,k_\sigma}F(\beta_\sigma^j)(x_\sigma^j).
\end{equation}
Consider the $(t+1)$-simplices for $j=1,..,k_\sigma$
$$ \sigma^j=\xymatrix{i_\sigma^j\ar[r]^{\beta_\sigma^j} & \sigma_0 \ar[r]^{\alpha_1}
& \sigma_1 \ar[r]^{\alpha_2}&...\ar[r]^{\alpha_{t-1}}&
\sigma_{t-1}\ar[r]^{\alpha_t} &\sigma_t}$$ and the $(t+1)$-chain
of $(M_p)^{i-1}_{t+1}$
$$ y_\sigma=\oplus_{j=1,..,k_\sigma}i_{\sigma^j}(x_\sigma^j).$$

The differential of $y_\sigma$ is
\begin{align}
d(y_\sigma)=&d_0(y_\sigma)+\sum_{j=1,..,t} (-1)^j d_j(y_\sigma)\notag\\
=&d_0(y_\sigma)+R_\sigma  \textit{, where $R_\sigma=\sum_{j=1,..,t} (-1)^j d_j(y_\sigma)$}\notag\\
=&\sum_{j=1,..,k_\sigma}i_{d_0(\sigma^j)}(F(\beta_\sigma^j)(x_\sigma^j))+R_\sigma
\notag\\
=&\sum_{j=1,..,k_\sigma}i_{\sigma}(F(\beta_\sigma^j)(x_\sigma^j))+R_\sigma\notag\\
=&i_{\sigma}(\sum_{j=1,..,k_\sigma}F(\beta_\sigma^j)(x_\sigma^j))+R_\sigma\notag\\
=&i_{\sigma}((x'_i)_\sigma) +R_\sigma\notag
\end{align}

where the last equality is due to (\ref{xenim}). Notice that
$R_\sigma$ lives in the subgroup $\bigoplus_{\sigma \in {N\P}_t,
deg(\sigma_0)=i-1, deg(\sigma_t)=p} F_{\sigma}\subseteq
(M_p)^{p'}_t$ of simplices beginning at degree $(i-1)$. Repeating
the same construction for each $\sigma \in{N\P}_t$ with
$deg(\sigma_0)=i$ and $d_0(\sigma)=\epsilon$ we obtain
$y_\epsilon=\sum_{\sigma} y_\sigma$ such that
$$
d(y_\epsilon)=\oplus_{\sigma \in{N\P}_t, deg(\sigma_0)=i,
d_0(\sigma)=\epsilon} (x'_i)_{\sigma}+R_\epsilon
$$
where $R_\epsilon$ lives in the subgroup $\bigoplus_{\sigma \in
{N\P}_t, deg(\sigma_0)=i-1, deg(\sigma_t)=p} F_{\sigma}\subseteq
(M_p)^{p'}_t$. Repeating the same argument for every $\epsilon\in
{N\P}_{t-1}$ with $deg(\epsilon_0)>i $ and $deg(\epsilon_{t-1})=p$
we obtain $y=\sum_\epsilon y_\epsilon$ such that
$$
d(y)=\oplus_{\sigma \in{N\P}_t, deg(\sigma_0)=i, d_0(\sigma_t)=p}
(x'_i)_{\sigma}+R
$$
where $R$ lives in the subgroup $\bigoplus_{\sigma \in {N\P}_t,
deg(\sigma_0)=i-1, deg(\sigma_t)=p} F_{\sigma}\subseteq
(M_p)^{p'}_t$. By construction $y\in (M_p)^{i-1}_{t+1}\subseteq
(M_p)^0_{t+1}$ and $d(y)\in (M_p)^{i-1}_t\subseteq
(M_p)^{p'}_{t+1}$. Thus $d(y)\in (M_p)^{p'}_t\cap
d((M_p)^0_{t+1})$. Then, by (\ref{E_p_r_big_dir}),
$[x'_i]=[x'_i-d(y)]=[x'_{i-1}]$ where $$x'_{i-1}=\oplus_{\sigma
\in{N\P}_t, i>deg(\sigma_0)\geq p', d_0(\sigma_t)=p}
(x'_i)_{\sigma}+ R$$ is a representative that lives in
$$\bigoplus_{\sigma \in {N\P}_t, i-1\geq deg(\sigma_0)\geq p',
deg(\sigma_t)=p} F_{\sigma}\subseteq (M_p)^{p'}_t$$ as wished.
That condition (\ref{condition_claim_step_2}) holds is clear from
the definition of $x'_{i-1}$.

\textbf{Step 3:} By Step $2$ we can suppose that
$$
x=\bigoplus_{\sigma \in {N\P}_t, deg(\sigma_0)=p' ,
deg(\sigma_t)=p} x_{\sigma}.$$ Our objective now is to see that
there exists $y\in (M_p)^0_{t+1}$ with $d(y)=x$. This implies that
$[x]=0$ and finishes the proof of the theorem. We need the
\begin{Claim}
\label{claim_step_3} There exist chains $x_i\in (M_p)^0_t$ for
$i=p',..,0$ and $y_i \in (M_p)^0_{t+1}$ for $i=p',..,1$ such that
\begin{equation}\label{equ_claim_step_3}
d(y_i)=x_i+x_{i-1}
\end{equation}
for $i=p',..,1$ with $x_{p'}=x$ and $x_0=0$ such that
\begin{enumerate} \item \label{claim_step_3_cond_1} $x_i$ lives on
$\bigoplus_{\sigma \in {N\P}_t, deg(\sigma_0)=i, deg(\sigma_t)=p}
F_{\sigma}\subseteq (M_p)^0_t$ for $i=p',..,0$. \item
\label{claim_step_3_cond_2} $d(x_i)=0$ for $i=p',..,0$.
\end{enumerate}
\end{Claim}
Notice that the claim finishes Step $3$: as $x_0=0$ then
$x_1=d(y_1)$, $x_2=d(y_2)-x_1=d(y_2-y_1)$,
$x_3=d(y_3)-x_2=d(y_3-y_2+y_1)$,..,
$x=x_{p'}=d(y_{p'})-x_{p'-1}=d(y_{p'}-y_{p'-1}+...+(-1)^{p'+1}y_1)$
where $y_{p'}-y_{p'-1}+...+(-1)^{p'+1}y_1\in (M_p)^0_{t+1}$.

Define $x_{p'}\definicio x$. Then condition
(\ref{claim_step_3_cond_1}) and (\ref{claim_step_3_cond_2}) are
satisfied for $i=p'$. We construct $y_i$ and $x_{i-1}$ from $x_i$
recursively beginning on $i=p'$. The arguments are similar to
those used in step $2$.

The differential $d$ on $L^p/L^{p-1}$ restricts to
$$d:(M_p)^0_t\rightarrow (M_p)^0_{t-1}.$$
As $d(x_{p'})=d(x)=0$, for every $\epsilon\in {N\P}_{t-1}$ with
$deg(\epsilon_{t-1})=p$ we can apply the projection
$$\pi_\epsilon: (M_p)^0_{t-1} \twoheadrightarrow F_\epsilon$$ and
obtain $\pi_\epsilon(d(x))=0$. If $deg(\epsilon_0)>p'$ then
$$\pi_\epsilon(d(x))=\sum_{\sigma \in{N\P}_t, d_0(\sigma)=\epsilon}
F(\alpha_1)(x_{\sigma})
$$
and thus
\begin{equation}\label{equ_epsilon_step_3}
0=\sum_{\sigma \in{N\P}_t, d_0(\sigma)=\epsilon}
F(\alpha_1)(x_{\sigma})
\end{equation}
for each $\epsilon\in {N\P}_{t-1}$ with $deg(\epsilon_0)>p'$ and
$deg(\epsilon_{t-1})=p$. Notice that each summand $x_\sigma$ with
$\sigma \in {N\P}_t$, $deg(\sigma_0)=p'$ and $deg(\sigma)=p$
appears in one and just one equation as (\ref{equ_epsilon_step_3})
(take $\epsilon=d_0(\sigma)$). Using now pseudo-injectivity we
build as before $y_\sigma$, $y_\epsilon=\sum_{\sigma} y_\sigma$
and $y=\sum_{\epsilon} y_\epsilon$, where $\epsilon$ runs over
$\epsilon\in {N\P}_{t-1}$ with $deg(\epsilon_0)>p'$ and
$deg(\epsilon_{t-1})=p$, such that
$$
d(y)=x+R
$$
with $R$ living in $\bigoplus_{\sigma \in {N\P}_t,
deg(\sigma_0)=p'-1, deg(\sigma_t)=p} F_{\sigma}\subseteq
(M_p)^0_t$. Call $y_{p'}\definicio y$ and $x_{p'-1}=R$. Then
Equation (\ref{equ_claim_step_3}) is satisfied. Condition
(\ref{claim_step_3_cond_1}) for $i=p'-1$ holds by the construction
of $R$ and condition (\ref{claim_step_3_cond_2}) for $i=p'-1$
holds because $d(x_{p'-1})=d(R)=d(d(y)-x)=d^2(y)-d(x)=0-0=0$ as
$d$ is a differential and $d(x)=0$ by hypothesis. The construction
of $y_i$ and $x_{i-1}$ from $x_i$ is totally analogous to the
construction of $y_{p'}$ and $x_{p'-1}$ from $x_{p'}$ that we have
just made.

After we have built $y_1$ and $x_0$ if we try to build
$y=\sum_\epsilon y_\epsilon$ and $R$ from $x_0$ we find that,
because there are not objects of negative degree (thus if $z\in
Im(i')$ where $deg(i')=0$ then $z=0$), $x_0=0$.
\end{proof}

The following examples come from Example \ref{examples_projective}.
They show the weaker conditions that are needed for
$\limdir$-acyclicity instead of projectiveness.

\begin{Ex}\label{limacyI_expushout}
For the ``pushout category" $\P$ with shape
$$\xymatrix {a_0 \ar[r]^{f}\ar[d]^{g} & b_1 \\ c_1}$$
a functor $F:\P\rightarrow \Ab$ is $\limdir$-acyclic if $F(f)$ and
$F(g)$ are monomorphisms.

For the ``telescope category" $\P$ with shape
$$ \xymatrix{a_0 \ar[r]^{f_1} & a_1 \ar[r]^{f_2} & a_2 \ar[r]^{f_3} &  a_3 \ar[r]^{f_4} & a_4 ...}$$
a functor $F:\P\rightarrow \Ab$ is $\limdir$-acyclic if
$F(f_{i}\circ f_{i-1} \circ ..\circ f_1)$ is a monomorphism and
$\ker F(f_{i}\circ f_{i-1} \circ ..\circ f_{i-d+1})\subseteq \im
F(f_{i-d})$ for $d=1,2,3,..,i-1$  for each $i=2,3,4,...$

Notice that for this it is enough that $F(f_i)$ is a monomorphism
for each $i=1,2,3,..$.
\end{Ex}

%% file: dual_results.tex
\section{Dual results for injective objects in $\Ab^\P$.}\label{section_dual_results}
The appropriate notions to characterize the injective objects in the
functor category $\Ab^\P$, where $\P$ is a graded poset, are the
following:

\begin{Defi}\label{defi_ker}
$\ker_F(i_0)=\bigcap_{i_0\stackrel{\alpha}\rightarrow i,\alpha\neq
1_{i_0}} \ker F(\alpha)$ (or $\ker_F(i_0)=F(i_0)$ if the index set
of the intersection is empty) and $\coim_F(i_0)=F(i_0)/\ker_F(i_0)$.
\end{Defi}

\begin{Defi}\label{property_pseudo-injective}
Let $F:\P\rightarrow \Ab$ be a functor over a graded poset $\P$ with
degree function $deg$. Fix an integer $d\geq 0$. If for any object
$i_0$ of $\P$, different objects $\{i_j\}_{j\in J}$ of $\P$, arrows
$\alpha_j:i_0\rightarrow i_j$ with $deg(\alpha_j)=d$ and elements
$x_j\in \ker_F(i_j)$ for each $j\in J$, there is $y\in F(i_0)$ with
$$F(\alpha_j)(y)=x_j$$
for each $j\in J$, we call $F$ \emph{$d$-pseudo-injective}. If $F$
is $d$-pseudo-injective for each $d\geq 0$ we call $F$
pseudo-injective.
\end{Defi}

Then we can prove the following

\begin{Thm}\label{pro_injective_bounded_graded}
Let $\P$ be a bounded above graded poset and $F:\P\rightarrow \Ab$
be a functor. Then $F$ is injective if and only if
\begin{enumerate}
\item for any object $i_0$ of $\P$ $\ker_F(i_0)$ is injective in $\Ab$.
\label{Defi_preinjective_1}\item $F$ is pseudo-injective.
\label{Defi_preinjective_2}
\end{enumerate}
\end{Thm}

Also in the dual case pseudo-injectiveness  is enough for vanishing
higher inverse limits:

\begin{Defi}\label{Defi_inverse_acyclic}
Let $\P$ be a graded poset and $F:\P\rightarrow \Ab$. We say $F$ is
\emph{$\liminv$-acyclic} if $\liminv^i F=0$ for $i\geq 1$.
\end{Defi}

\begin{Thm}\label{pro_inv_acyclic_graded}
Let $F:\P\rightarrow \Ab$ be a pseudo-injective functor over a
bounded above graded poset $\P$. Then $F$ is $\liminv$-acyclic.
\end{Thm}